\documentclass[12pt]{article}
\usepackage{amssymb}
\usepackage{graphicx}
\usepackage{amsfonts}
\usepackage{latexsym}
\usepackage{amsthm}
\usepackage{amsmath}
\usepackage[all,arc]{xy} 
\usepackage{amsmath, amssymb, graphics, setspace}

\usepackage{amssymb, amscd}

\usepackage{url}

\newtheorem{problem}{Problem}
\newtheorem{theo}[problem]{Theorem}

\newtheorem{defin}[problem]{Definition}
\newtheorem{prop}[problem]{Proposition}

\newtheorem{lemma}[problem]{Lemma}
\newtheorem{exam}[problem]{Example}

\begin{document}
\date{September 9, 2014}
 \title{{Signed polyomino tilings by \\ $n$-in-line polyominoes and Gr\" obner bases}}

\author{{Manuela Muzika Dizdarevi\' c}\\ {\small Faculty of Natural Sciences}\\[-2mm] {\small and Mathematics, Sarajevo}
\and Marinko Timotijevi\'{c} \\ {\small Department of Mathematics and Informatics, }\\[-2mm]{\small Faculty of Science, University of Kragujevac }
\and Rade  T.\ \v Zivaljevi\' c\\ {\small Mathematical Institute}\\[-2mm] {\small SASA, Belgrade} }

\maketitle
\begin{abstract}
Conway and Lagarias observed that a triangular region $T(m)$ in a
hexagonal lattice admits a {\em signed tiling} by three-in-line
polyominoes (tribones) if and only if $m\in \{9d-1, 9d\}_{d\in
\mathbb{N}}$. We apply the theory of Gr\"{o}bner bases over
integers to show that $T(m)$ admits a signed tiling by $n$-in-line
polyominoes ($n$-bones) if and only if $$m\in \{dn^2-1,
dn^2\}_{d\in \mathbb{N}}.$$ Explicit description of the
Gr\"{o}bner basis allows us to calculate the `Gr\"{o}bner discrete
volume' of a lattice region by applying the division algorithm to
its `Newton polynomial'. Among immediate consequences is a
description of the {\em tile homology group} of the $n$-in-line
polyomino.

\end{abstract}

\renewcommand{\thefootnote}{\fnsymbol{footnote}}
\footnotetext{R.~\v Zivaljevi\'c was supported by the Grants
174017 and 174020 of the Ministry for Science and Technological
Development of Serbia.}


\section{Introduction}\label{sec:main}

A $n$-bone is by definition a $n$-in-line polyomino (polyhex) in a
hexagonal lattice. For example a $3$-bone is the same as the {\em
tribone} in the sense of \cite{Thur}. One initial objective is to
determine when a triangular region $T(m)$ in a hexagonal lattice
admits a signed tiling by $n$-bones.

\begin{figure}[hbt]
\centering
\includegraphics[scale=0.60]{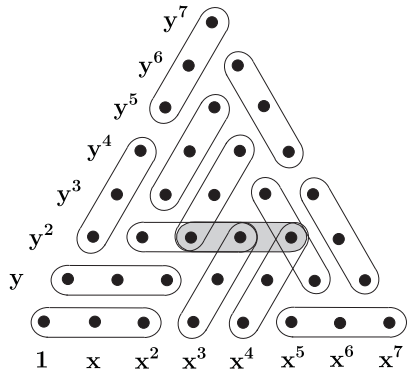}
\caption{A signed tiling of a triangular region by $3$-bones.}
\label{fig:trimino-5}
\end{figure}

\medskip
By a theorem of Conway and Lagarias (\cite[Theorem~1.2.]{ConLag})
$T(m)$ admits a signed tiling by $3$-bones if and only if $m=9d$
or $m=9d-1$ for some integer $d\geq 1$, the case $m=8$ is
exhibited in Figure~\ref{fig:trimino-5}. Our central result is
Theorem~\ref{thm:main} which claims that $T(m)$ admits a signed
tiling by $n$-bones if and only if $m=dn^2$ or $m=dn^2-1$ for some
integer $d\geq 1$.

\medskip
The Gr\"{o}bner basis approach to signed polyomino tilings was
originally proposed by Bodini and Nouvel \cite{BN}, see also
\cite{MDZ} for an application to tilings with symmetries. The
knowledge of the Gr\"{obner} basis (Theorem~\ref{thm:G-main})
offers a deeper insight into the (signed) tiling problem  and
provides a powerful tool for analyzing general behavior and
selected particular cases. It is well adopted to other methods of
lattice geometry and we illustrate this by examples involving
Brion's theorem (Example~\ref{exam:Brion}).

\medskip
Computing Gr\"{o}bner basis of a tiling problem yields as a
byproduct complete information about the associated {\em tile
homology group} \cite{ConLag,Reid}. In general computing homology
class by a `division algorithm' may offer an interesting new
computational paradigm which deserves further exploration.

\section{Gr\"{o}bner bases}

The notion of a {\em strong Gr\"{obner} base} \cite{AL, Licht}
(called a $D$-Gr\"{o}bner base in \cite{BW}) allows us to apply
the Gr\"{o}bner basis theory to polynomials with integer
coefficients. Here is a brief outline of some basic definitions
and theorems with pointers to some of the key references.

\medskip
A term is a product $t=cx^\alpha$ where $c$ is the coefficient and
$x^\alpha = x_1^{\alpha_1}\cdots x_k^{\alpha_k}$ is the associated
monomial (power product). For a given polynomial $f\in
\mathbb{Z}[x_1,x_2,\ldots, x_k]$ the associated remainder on
division by a Gr\"{o}bner basis $G$ is $\overline{f}^G$ and $f$
reduces to zero $f \stackrel{G}{\longrightarrow} 0$ if
$\overline{f}^G = 0$. $LM(f)$ and $LC(f)$ are respectively the
leading monomial and the leading coefficient with respect to the
chosen term order $\preceq$. We write $lcm(a,b)$ and $gcd(a,b)$
respectively for the least common multiple and the greatest common
divisor of $a$ and $b$.

For other basic notions of Gr\"{o}bner basis theory (over
integers), such as $S$-polynomial, standard representation etc.,
the reader is referred to \cite{AL, BW, Licht} (see also
\cite{CLO1, CLO2, Sturm} for related results for coefficients in a
field).

\subsection{Gr\"{o}bner bases over principal ideal domains}

Let $\Lambda = R[x_1,\ldots, x_k]$ be the ring of polynomials with
coefficients in a principal ideal domain $R$. For a given ideal
$I\subset\Lambda$ the associated {\em strong Gr\"{o}bner basis},
called also the $D$ bases in \cite{BW}, may be introduced as
follows (see \cite[p.~251]{AL} and \cite[p.~455]{BW}).

\begin{defin}\label{def:grobner}
A finite set $G\subset I$ is a strong Gr\"{o}bner basis of $I$
(with respect to the chosen term order $\preceq$) if for each
$f\in  I\setminus\{0\}$ there exists $g\in G$ such that the
leading term of $f$ is divisible by the leading term of $g$,
$LT(g) | LT(f)$, meaning  that $LT(f) = t  LT(g)$ for some term
$t$.
\end{defin}

The following theorem provides a useful criterion for testing
whether a finite set of polynomials is a Gr\"{o}bner basis of the
ideal generated by them, see \cite[Chapter~10,
Corollary~10.12]{BW} and \cite[Theorem~2.1.]{bane}.

\begin{theo}\label{thm:criterion}
Let $G$ be a finite collection of non-zero polynomials which
generate an ideal $I_G$. Suppose that,
 \begin{enumerate}
 \item[{\em (1)}] For each pair $g_1, g_2\in G$ there exists $h\in
 G$ such that, $$LM(h)\vert {\rm lcm}(LM(g_1), LM(g_2)) \mbox{ {\rm and} } LC(h)\vert {\rm gcd}(LC(g_1),
 LC(g_2))$$
 \item[{\em (2)}] For each pair $g_1, g_2\in G$ the associated $S$-polynomial reduces to zero,
 $$S(g_1,g_2)\stackrel{G}{\longrightarrow} 0.$$
 \end{enumerate}
Then $G$ is a strong Gr\"{o}bner basis of $I_G$.
\end{theo}

\subsection{Gr\"{o}bner bases over Euclidean domains}

The general theory is further simplified if one works with
Euclidean domains. Aside from standard references \cite{AL, BW} a
self-contained account can be found in \cite{Licht}.  In the case
of integers one usually chooses the linear ordering,
\begin{equation}
\label{eqn:ordering-1} \ldots   0 < +1 < -1 < +2 < -2 < +3 < - 3 <
\ldots
\end{equation}
which allows us to define unambiguously remainders,
$S$-polynomials etc. For example following (\ref{eqn:ordering-1})
the reduction of $8$ mod $5$ is $-2$ rather than $+3$.

\medskip\noindent
{\bf Caveat:} We find it convenient in
Section~\ref{sec:remainders} to stick to positive remainders and
write that $+3$ is, rather than $-2$, the remainder of $8$ on
division by $5$. In other words we use the following term order
for coefficients,
\begin{equation}
\label{eqn:ordering-2} \ldots   0 < +1 < +2 < +3 <  \ldots < -1 <
-2 < - 3 < \ldots .
\end{equation}

\begin{exam}\label{ex:1} {\rm
In agreement with (\ref{eqn:ordering-1}) many standard computer
algebra packages (including Wolfram Mathematica 9.0) would yield
$-1 -x - y$ as the remainder of $T(6)$
(Section~\ref{sec:n-bone-ideal}) on division by $GBI_3$. In
Section~\ref{sec:remainders} we would (following
(\ref{eqn:ordering-2})) reduce this polynomial further by the
element $g_3(3)=3T(2)$ (Section~\ref{sec:GB-n-bone}) and obtain
the polynomial $2+x+y$. }
\end{exam}

\section{From polyominoes to polynomials}

Each polyomino $P\subset \mathbb{Z}^2$ is associated the
corresponding `Newton polynomial' $f_P := \sum_{(p,q)\in
P}~x^py^q$. For example the shaded tribone $P$ in
Figure~\ref{fig:trimino-5} is associated the trinomial
$x^2y^2+x^3y^2+x^4y^2$.

\begin{prop}\label{prop:poly-tiling}
A polyomino $P$ admits a {\em signed tiling} by translates of {\em
prototiles} $P_1, P_2, \ldots, P_k$ if and only if for some (test)
monomial $x^\alpha = x_1^{\alpha_1}\ldots x_n^{\alpha_n}$ the
polynomial $x^\alpha f_P$ is in the ideal generated by polynomials
$f_{P_1}, \ldots , f_{P_k}$,
\begin{equation}\label{prop:criterion-for-tiling}
x^\alpha f_P \in \langle f_{P_1}, f_{P_2} \ldots , f_{P_k}
\rangle.
\end{equation}
Moreover, the set of test monomials $\mathcal{T} = \{x^\alpha \mid
\alpha\in T\}$ can be chosen from any set $T\subset \mathbb{N}^n$
of multi-indices which is {\em cofinal} in $(\mathbb{N}, \leq)$.
\end{prop}

\medskip\noindent
{\bf Proof:} Let $J\subset \mathbb{Z}[x,y; x^{-1}, y^{-1}]$ be the
extension of the ideal $\langle f_{P_1}, f_{P_2} \ldots , f_{P_k}
\rangle$  in the ring of Laurent polynomials with coefficients in
$\mathbb{Z}$. $P$ admits a signed tiling by translates of
prototiles $P_1, P_2, \ldots, P_k$ if and only if $f_P\in J$. The
proposition is an immediate consequence of the relation,
\[
  J = \bigcup_{x^\alpha\in \mathcal{T}} x^{-\alpha} \langle f_{P_1}, f_{P_2} \ldots , f_{P_k}
\rangle.
\]

\section{The $n$-bone ideal $I_n$}
\label{sec:n-bone-ideal}

Let $I_n = \langle b_1(n), b_2(n), b_3(n) \rangle \subset
\mathbb{Z}[x,y]$ be the ideal generated by polynomials,
\begin{equation}\label{eqn:B-polis}
b_1(n) = 1+x+\ldots + x^{n-1},  b_2(n) = 1+y+\ldots + y^{n-1},
b_3(n) = x^{n-1}+x^{n-2}y+\ldots + y^{n-1}
\end{equation}
These polynomials correspond to three types of $n$-in-line
polyominoes in a hexagonal lattice.

\medskip
We denote by $T(m)$ the `integer-point transform' \cite[p.~60]{BR}
(Newton polynomial) of a triangular region with the side-length
equal to $m$,
\begin{equation} T(m) =
\sum\limits_{0\leq i,j\leq m-1 \atop i+j\leq m-1}x^{i}y^{j}.
\end{equation}

\section{Gr\"{o}bner basis for the $n$-bone ideal}
\label{sec:GB-n-bone}

Let $GBI_n = \{g_1(n), g_2(n), g_3(n), g_4(n)\}$ be the following
set of polynomials,
\begin{equation}
\begin{array}{ccl}
g_1(n) & = & b_1(n)\\
g_2(n) & = & b_2(n)\\
g_3(n) & = & nT(n-1)\\
g_4(n) & =&b_3(n)-b_1(n)-b_2(n)
\end{array}
\end{equation}
\begin{lemma}\label{lemma:LT}
The leading terms of polynomials $g_1, g_2, g_3, g_4$ with respect
to the lexicographical term order are the following,
\begin{equation}\label{eqn:LT}
LT(g_1(n))= x^{n-1},\, LT(g_2(n))= y^{n-1},\, LT(g_3(n))=
nx^{n-2}, \, LT(g_4(n))= x^{n-2}y
\end{equation}
\end{lemma}
The relations listed in Proposition~\ref{prop:1} will be needed in
the sequel. The first equality is trivial while the rest follow
from an iterated application of the identity $a^d-b^d =
a^{d-1}+a^{d-2}b+\ldots+b^{d-1}$ for suitable $a$ and $b$.
\begin{prop}\label{prop:1}
\begin{equation*}
\begin{array}{rcl}
T(n)&=& T(n-1)+b_3(n)\\
(x-1)T(n-1)&=&b_3(n)-b_2(n)\\
(x-y)T(n-1)&=&b_1(n)-b_2(n)\\
(y-1)g_1(n) + (y-x)g_4(n)&=&(x-1)g_2(n).
\end{array}
\end{equation*}
\end{prop}

\begin{prop}\label{prop:2}
The set $GBI_n$ is a basis of the ideal $I_n$.
\end{prop}

\medskip\noindent
{\bf Proof:} Let $\langle GBI_n\rangle$ be the ideal generated by
$GBI_n$. It is obvious that
$$
I_n = \langle g_1(n), g_2(n), g_4(n) \rangle \subseteq \langle
GBI_n \rangle
$$
so it is sufficient to show that $g_3(n)\in I_n$. As a consequence
of the second identity in Proposition~\ref{prop:1},
\begin{equation*}
\begin{array}{rcl}
(x-1)T(n-1)&\in& I_n\\
(x^2-1)T(n-1)&\in& I_n\\
&\vdots&\\
(x^{n-1}-1)T(n-1) &\in& I_n
\end{array}
\end{equation*}
By adding these polynomials we obtain
\begin{equation*}
b_1(n)T(n-1)-nT(n-1) \in I_n
\end{equation*}
and $g_3=nT(n-1)\in I_n$ which is the desired conclusion. \hfill
$\square$

\begin{theo}\label{thm:G-main}
The set of polynomials $GBI_n$ is a strong Gr\"{o}bner basis (over
the base ring $ \mathbb{Z}$) of the ideal $I_n$, $n \geq 2$, with
respect to lexicographic term order.
\end{theo}

\medskip\noindent
{\bf Proof:} The case $n=2$ is elementary so we assume that $n\geq
3$. By Proposition~\ref{prop:2} the set $GBI_n$ is a basis of the
ideal $I_n$. In order to show that this is indeed a strong
Gr\"{o}bner basis of the ideal $I_n\subset\mathbb{Z}[x,y]$ we
apply the $\mathbb{Z}$-version of the Buchberger criterion.

\medskip
Following  \cite[Theorem~2]{Licht} it is sufficient to show that
for every pair of polynomials $g_i(n), g_j(n)\in GBI_n$, their
$S$-polynomial reduces to $0$ by the set $GBI_n$. Equivalently,
one can use Theorem~\ref{thm:criterion} by observing that the
condition (1) is (in Light of Lemma~\ref{lemma:LT}) readily
satisfied.

\medskip
Since the leading monomials of polynomials $g_1(n), g_2(n)$ and
$g_2(n), g_3(n)$ are pairwise coprime (Lemma~\ref{lemma:LT}) and
the leading coefficients divide each other, from \cite[Theorem
3]{Licht} we conclude that
\begin{equation*}
S(g_1(n),g_2(n))\stackrel{GBI_n}{\xrightarrow{\hspace*{1.2cm}}} 0
\quad\mbox{and} \quad S(g_2(n),g_3(n))
\stackrel{GBI_n}{\xrightarrow{\hspace*{1.2cm}}} 0 .
\end{equation*}
Let us consider now polynomials $g_1(n)$ and $g_4(n)$. By
Lemma~\ref{lemma:LT},
\begin{equation*}
S(g_1(n),g_4(n))=yg_1(n)-xg_4(n).
\end{equation*}
Since
\begin{equation*}
LT(S(g_1(n),g_4(n)))=LT(x^{n-1}+x^{n-2}y-x^{n-2}+\dots)=x^{n-1}
\end{equation*}
we can reduce this polynomial by $g_1(n)$. The reduction leads to
the polynomial,
\begin{equation*}
S(g_1(n),g_4(n))-g_1(n)=yg_1(n)-xg_4(n) - g_1(n)
\end{equation*}
which has the leading term
\begin{equation*}
LT(S(g_1(n),g_4(n))-g_1(n))=LT(-x^{n-2}y^2+x^{n-2}y-\dots)=-x^{n-2}y^{2}
\end{equation*}
and which, in light of Lemma~\ref{lemma:LT}, can be reduced by
$g_4(n)$. This reduction leads to the polynomial,
\begin{equation*}
S(g_1(n),g_4(n))-g_1(n)+yg_4(n)=(y-1)g_1(n)+(y-x)g_4(n).
\end{equation*}
By using the last equality in Proposition \ref {prop:1} we finally
get a strong representation of $S(g_1(n),g_4(n))$ by the set
$GBI_n$,
\begin{equation}\label{eqn:case14}
S(g_1(n),g_4(n))=g_1(n)-yg_4(n)+(x-1)g_2(n).
\end{equation}
In a similar manner we show reducibility of polynomials
$S(g_1(n),g_3(n))$ and $S(g_3(n),g_4(n))$.

\medskip
By Lemma~\ref{lemma:LT}, $S(g_1(n), g_3(n))= ng_1(n)-xg_3(n)$ has
the leading term $-nx^{n-2}y$. Consequently it can be reduced by
the polynomial $g_4(n)$ and we focus our attention to the
polynomial,
\[
ng_1(n) - xg_3(n) + ng_4(n).
\]
This polynomial is reducible to zero since, in light of the second
equality in Proposition~\ref{prop:1}, it is equal to $-ng_3(n)$.
In particular it has the strong representation in terms of the
basis $GBI_n$,
\[
S(g_1(n),g_3(n)) =  -ng_4(n)-g_3(n).
\]
A similar calculation shows that,
\[
S(g_3(n),g_4(n)) = g_3(n)+ng_2(n)
\]
is a strong representation of $S(g_3(n),g_4(n))$.
\medskip
Together with the case of $S$-polynomial $S(g_2(n),g_4(n))$, which
is separately treated in Lemma~\ref{lemma:Sg2g4}, this concludes
the proof of Theorem~\ref{thm:G-main}. \hfill $\square$

\begin{figure}[hbt]
\centering
\includegraphics[scale=0.50]{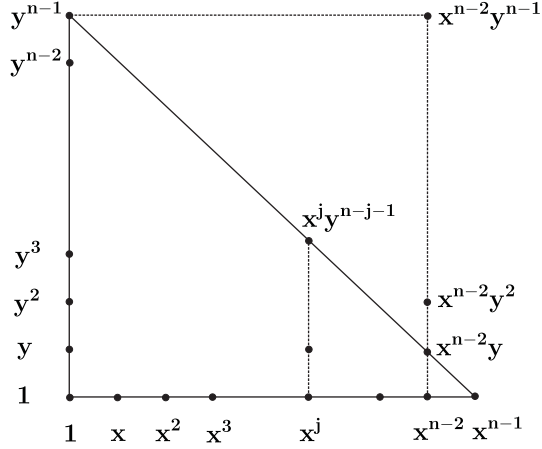}
\caption{Reduction of $S(g_2(n), g_4(n))$.} \label{fig:S-poli-1}
\end{figure}

\begin{lemma}\label{lemma:Sg2g4}
The $S$-polynomial $S(g_2(n), g_4(n))$ can be reduced to $0$ by
the basis $GBI_n$.
\end{lemma}

\medskip\noindent
{\bf Proof:} By Lemma~\ref{lemma:LT}, $S(g_2(n), g_4(n))=
x^{n-2}g_2(n) - y^{n-2}g_4(n)$. The terms $x^{n-2}y$ and
$-x^{n-2}$ are the leading two terms of the polynomial $g_4(n)$
and they are the only terms in the lexicographically leading
column $\{x^{n-2}y^i\}_{i\geq 0}$ (Figure~\ref{fig:S-poli-1}).
This observation indicates that one should begin with the
reduction of the $S$-polynomial $S(g_2(n), x^{n-2}y-x^{n-2}) =
x^{n-2}S(b_2(n),y-1)$. From the identity,
\begin{equation}\label{eqn:again}
b_2(n) - n = \sum_{j=0}^{n-1}~(y^j-1) = (y-1)B_2(n)
\end{equation}
where $B_2(n) = b_2(n-1)+b_2(n-2)+\ldots + b_2(1)$ we observe that
$S(g_2(n), g_4(n))$ can be reduced to the polynomial
$x^{n-2}g_2(n) - B_2(n)g_4$ which has the monomial $nx^{n-2}$ as
the leading term. This is precisely the leading term of the
polynomial $g_3(n) = nT(n-1)$ so we turn our attention to the
polynomial,
\begin{equation}\label{eqn:attention}
x^{n-2}g_2(n) - B_2(n)g_4(n) - g_3(n)
\end{equation}
Since by definition $b_3(n)-b_1(n) =
\sum_{k=1}^{n-1}~x^{n-k-1}(y^k-1)$ we observe (in light of
(\ref{eqn:again})) that,
\[
B_2(n)[b_3(n)-b_1(n)] = [\sum_{k=1}^{n-1}
x^{n-k-1}(\sum_{j=0}^{k-1}~y^j)][b_2(n)-n] = T(n-1)b_2(n) -
nT(n-1).
\]
It follows that
$$
B_2(n)g_4 + g_3(n) = [T(n-1) - B_2(n)]b_2(n)
$$
which implies that the polynomial (\ref{eqn:attention}) can be
reduced by $g_2(n) = b_2(n)$ with zero remainder. This completes
the proof of the lemma. \hfill $\square$

\section{Evaluation of remainders }
\label{sec:remainders}

Our objective in this section is to calculate the reminder
$\overline{T(n)}^{GBI_n}$ of $T(n)$ on division by the Gr\"{o}bner
basis $GBI_n$.

\begin{lemma}\label{lemma:division}
Suppose that
\begin{equation}\label{eqn:deljenje-po-x}
p(x) = q(x)(x^n-1)+r(x)
\end{equation}
is the equality arising from the division of a polynomial $p(x)\in
\mathbb{Z}[x]$ by $x^n-1$ where $q(x)$ is the quotient and $r(x)$
the remainder.

If $P(x,y) = \frac{p(x)-p(y)}{x-y}$ and $R(x,y) =
\frac{r(x)-r(y)}{x-y}$ then,
\begin{equation}
\overline{P(x,y)}^{GBI_n} = \overline{R(x,y)}^{GBI_n}.
\end{equation}
Moreover, if $R(x,y)$ cannot be further reduced by the Gr\"{o}bner
basis $GBI_n$ then the remainder of $P(x,y)$ on division by
$GBI_n$ is,
\begin{equation}
\overline{P(x,y)}^{GBI_n} = \overline{R(x,y)}^{GBI_n} = R(x,y) =
\frac{r(x)-r(y)}{x-y}.
\end{equation}
\end{lemma}

\medskip\noindent
{\bf Proof:} From (\ref{eqn:deljenje-po-x}) we deduce the
following equality,
\begin{equation}
\frac{p(x)-p(y)}{x-y} =  \frac{q(x)-q(y)}{x-y}(x^n-1) +
q(y)\frac{x^n-y^n}{x-y} + \frac{r(x)-r(y)}{x-y}.
\end{equation}
Both $x^n-1 = (x-1)b_1(n)$ and $\frac{x^n-y^n}{x-y} = b_3(n)$ are
in the ideal $I_n$ so $\overline{P(x,y)}^{GBI_n} =
\overline{R(x,y)}^{GBI_n}$. The second part of the lemma is an
immediate consequence. \hfill $\square$

\begin{lemma}\label{lemma:b3-ostatak}
Let $b_3(m)=x^{m-1}+x^{m-2}y+\ldots + y^{m-1}$ and assume by
convention that $b_3(0)=0$. Then,
\begin{equation}\label{eqn:ostatak}
\overline{b_3(m)}^{GBI_n} = b_3(r)
\end{equation}
where $r = r_m^n = m - \lfloor m/n \rfloor n$ is the reminder of
the division of $m$ by $n$.
\end{lemma}

\medskip\noindent
{\bf Proof:} Observe that $b_3(m) = P(x,y) =
\frac{p(x)-p(y)}{x-y}$ for $p(x) = x^{m}$. For this choice of
$p(x)$ the equation corresponding to (\ref{eqn:deljenje-po-x}) is
\[
x^m = (x^{m-n} + x^{m-2n}+\ldots + x^r)(x^n-1) + x^r.
\]
Since $LT(R(x,y))= LT(b_3(r))= x^{r-1}$ is not divisible by any of
the leading terms of the Gr\"{o}bner basis $GBI_n$ listed in
(\ref{eqn:LT}) we observe that $\overline{b_3(r)}^{GBI_n} =
b_3(r)$ and the result follows from the second half of
Lemma~\ref{lemma:division}. \hfill $\square$

\medskip
Since,
\begin{equation}\label{eqn-rekurent}
T(m) = T(m-1)+b_3(m)
\end{equation}
Lemma~3 may be used for an inductive evaluation of
$\overline{T(m)}^{GBI_n}$. As before $r = r_m = r_m^n = m-\lfloor
m/n \rfloor n$.

\begin{prop}\label{prop:ostatak-T(n)}
For each integer $n\geq 1$ the sequence of polynomials $\alpha_m^n
= \alpha_m^n(x,y) = \overline{T(m)}^{GBI_n}$ is periodic with the
period $n^2$.

For $1\leq m\leq n^2-2$, $T(m) = \sum_{k=1}^m~b_3(k)$ and
\begin{equation}\label{eqn:suma-T(m)}
\overline{T(m)}^{GBI_N} = \sum_{k=1}^m~b_3(r_{k}^n)\neq 0.
\end{equation}
For $m\in \{n^2-1, n^2\}$,
\begin{equation}\label{eqn:suma-T(m)-jos}
\overline{T(m)}^{GBI_N} =  0.
\end{equation}
\end{prop}

\medskip\noindent
{\bf Proof:} To establish the periodicity of the sequence
$\alpha_m=\alpha_m^n =  \overline{T(m)}^{GBI_n}$ it is sufficient
to establish the equalities (\ref{eqn:suma-T(m)}) and
(\ref{eqn:suma-T(m)-jos}).

Indeed, assume that (\ref{eqn:suma-T(m)}) and
(\ref{eqn:suma-T(m)-jos}) are true and that $\alpha_m$ is periodic
with the period $n^2$ in the interval $[1, jn^2]$ for some integer
$j\geq 1$. For each $d\in [jn^2+1, (j+1)n^2]$,
\[
\alpha_d = \overline{T(d)}^{GBI_n} = \overline{A + B}^{GBI_n}
\]
where $A = T(jn^2)$ and $B = \sum_{k=jn^2+1}^d~b_3(k)$. Since by
the inductive hypothesis $\overline{A}^{GBI_n}=0$ we observe that
\[
\alpha_d = \overline{B}^{GBI_n} = \sum_{k=jn^2+1}^d~b_3(r_k) =
\sum_{k=1}^{d'}~b_3(r_k)
\]
where $d'= d - \lfloor d/n^2\rfloor n^2$ which proves that the
sequence $\alpha_m$ repeats the same pattern in the interval
$[jn^2+1, (j+1)n^2]$.

\medskip
Since $T(m) = \sum_{k=1}^m~b_3(k)$, in light of the equality
(\ref{eqn:ostatak}) it is not surprising that,
\[
\alpha_m = \overline{T(m)}^{GBI_n} =
\overline{\sum_{k=1}^m~b_3(r_{k}^n)}^{GBI_n}.
\]
The equality (\ref{eqn:suma-T(m)}) claims more than that, it says
that the right hand side ${\rm rhs}$-(\ref{eqn:suma-T(m)}) of
(\ref{eqn:suma-T(m)}) is reduced with respect to the Gr\"{o}bner
basis $GBI_n$. Indeed, for $m\leq n^2-2$ if $Cx^py^q$ is the
leading term of ${\rm rhs}$-(\ref{eqn:suma-T(m)}) then either
$p<n-2$ or $C\leq n-1$.

\medskip
A similar analysis shows that $\overline{T(n^2-1)}^{GBI_n} =
nT(n-1)=g_4(n)\in I_n$. This together with the fact $b_3(n^2)\in
I_n$ establish the equality (\ref{eqn:suma-T(m)-jos}). \hfill
$\square$

\section{Signed tilings by $n$-bones}  

\begin{theo}\label{thm:main} A triangular region $T(m)$ in a
hexagonal lattice admits a {\em signed tiling} by $n$-in-line
polyominoes ($n$-bones) if and only if
\begin{equation}\label{eqn:main}
m \equiv -1 \quad {\rm mod}\, n^2 \qquad {\rm or} \qquad   m
\equiv 0 \quad {\rm mod}\, n^2.
\end{equation}
\end{theo}

\medskip\noindent
{\bf Proof:} By Proposition~\ref{prop:criterion-for-tiling} it is
sufficient to check if at least one of the polynomials,
\[
T(m), \quad x^ny^nT(m), \quad x^{2n}y^{2n}T(m), \quad
x^{3n}y^{3n}T(m), \quad \ldots
\]
is in the ideal $I_n$ generated by $n$-bones. Since
$x^{kn}y^{kn}-1\in I_n$ for each $k$, the triangular region $T(m)$
admits a signed tiling by $n$-in-line polyominoes if and only if
$T(m)\in I_n$.

\medskip
By Proposition~\ref{prop:ostatak-T(n)} this happens if and only if
the condition (\ref{eqn:main}) is satisfied. This observation
completes the proof of the theorem. \hfill $\square$

\section{Tile homology groups and Brion's theorem}\label{sec:tile-homology}

For {\em tile homology groups} the reader is referred to
\cite{ConLag} and \cite{Reid}. The following result illustrates
how one can read off the tile homology group from the Gr\"{o}bner
basis.

\begin{prop}
The tile homology group of a polyomino with prototiles
$\mathcal{P}$ and the associated ideal $I = I_{\mathcal{P}}\subset
\mathbb{Z}[x_1,\ldots, x_k] = \mathbb{Z}[\overline{x}]$ can be
computed as the direct limit ${\rm colim}_{\alpha\in \mathbb{N}^k}
\mathcal{D}_\alpha$ where $\mathcal{D}_\alpha =
\mathbb{Z}[\overline{x}]/I$ and for $\alpha\leq\beta$, the
connecting map $\mathcal{D}_\alpha \stackrel{\times x^{\beta -
\alpha}}{\longrightarrow} \mathcal{D}_\beta$ is multiplication by
$x^{\beta - \alpha}$.
\end{prop}
It is clear that this direct system can be in principle calculated
if a Gr\"{o}bner basis of the ideal $I$ is known. In favorable
cases, such as the case of the $n$-in-line polyomino, all
connecting maps are isomorphism (see the proof of
Theorem~\ref{thm:main}). The following proposition is a direct
consequence of Lemma~\ref{lemma:LT} and the fact that
$\mathbb{Z}[x,y]/I$ is generated by monomials which are reduced
with respect to the Gr\"{o}bner basis.

\begin{prop}\label{prop:tile-homology}
The tile homology group of the $n$-in-line polyomino is isomorphic
to the group, $$\mathbb{Z}^{(n-1)(n-2)}\oplus
\mathbb{Z}/n\mathbb{Z}.$$
\end{prop}

\medskip\noindent
The knowledge of a short Gr\"{o}bner basis provides powerful
experimental tool which is particularly well adopted to methods of
lattice geometry. Theorem~\ref{thm:main} was discovered by
experiments which involved Brion's theorem. Indeed, Brion's
theorem and its relatives provide a short rational form for the
integer-point transform which is an ideal input for a division
algorithm. The following example from {\em Mathematica 9.0}
exhibits the short rational form for the Newton polynomial
(integer-point transform) of the triangular region $T(n)$.

\begin{exam}\label{exam:Brion}
\begin{doublespace}
\noindent\(\pmb{T[\text{n$\_$}]\text{:=}
\text{Together}\left[\frac{1}{(1-x)*(1-y)}+\frac{x{}^{\wedge}(n+1)}{(x-1)*(x-y)}+
\frac{y{}^{\wedge}(n+1)}{(y-1)*(y-x)}\right]}\)
\end{doublespace}
\end{exam}

\section{Gr\"{o}bner discrete volume}

Let $Q$ be a convex polytope with vertices in $\mathbb{N}^d$ and
let $f_Q$ be its Newton polynomial (integer-point transform). The
usual `discrete volume' of $Q$, defined in \cite{Barv, BR} as the
number of integer points inside $Q$, can be evaluated as the
remainder of $f_Q$ on division by the ideal
\[
I = \langle x_1-1, x_2-1, \ldots , x_d-1 \rangle.
\]
Let $J\subset \mathbb{Z}[x_1,\ldots, x_d]$ be an ideal, say the
ideal associated to a set $\mathcal{R}$ of prototiles in
$\mathbb{N}^d$. Let $G=G_J$ be the Gr\"{o}bner basis of $J$ with
respect to some term order. It may be tempting to ask (at least
for some carefully chosen ideals $J$) what is the geometric and
combinatorial significance of the remainder $\overline{f}_Q^G$ of
the integer-point transform $f_Q$ on division by the Gr\"{o}bner
basis $G$.

\begin{defin}\label{def:G-discr-vol}
The polynomial valued function $Q \mapsto \overline{f}_Q^G $ is
referred to as Gr\"{obner} or $G$-discrete volume of $Q$ with
respect to the Gr\"{o}bner basis $G$,
\end{defin}
The Definition~\ref{def:G-discr-vol} may look somewhat artificial
at first sight. Note however that the basic geometric idea of a
volume of a geometric object $Q$ involves approximation, or rather
exhaustion (tiling!) of $Q$ by a set of prototiles $\mathcal{R}$.
The fact that the $G$-volume is a polynomial valued (rather than
integer valued) function reflects  the idea that there may be more
than one object in $\mathcal{R}$ used for `measurements' of $Q$.

\medskip
As in the case of integer-point enumeration in polyhedra, Brion's
theorem is a powerful tool for calculation of the $G$-discrete
volume.  It may be expected that some aspects of Ehrhart theory
can be extended in an interesting way to Gr\"{o}bner volumes, in
particular the results from Section~\ref{sec:remainders} can be
interpreted as the evaluation of the $GBI_n$-discrete volume of
the triangular region $T(m)$.

\bigskip\noindent
{\bf Acknowledgements:} The symbolic algebra computations in the
paper were performed with the aid of {Wolfram Mathematica $9.0$}.

\end{document}